\documentclass[uslettersize, 10pt, conference]{ieeeconf}      

\IEEEoverridecommandlockouts                              
\usepackage{graphics} 
\usepackage{epsfig} 
\usepackage{mathptmx} 
\usepackage{times} 
\usepackage{amsmath} 
\usepackage{amssymb}  
\usepackage{amscd}

\usepackage[noadjust]{cite}

\newtheorem{theorem}{Theorem}
\newtheorem{lem}[theorem]{Lemma}

\newtheorem{fact}[theorem]{Fact}

\title{\LARGE \bf
Adjacency Criterion For Gradient Flow With Multiple Local Maxima}

\author{Xudong Chen$^{1}$ 
\thanks{$^{1}$X. Chen is with the ECEE department, University of Colorado at Boulder.
        {\tt\small xudong.chen@colorado.edu}}%
}

\begin{document}

\maketitle
\thispagestyle{empty}
\pagestyle{empty}

\begin{abstract}
In this paper, we investigate the geometry of a general class of gradient flows with multiple local maxima. we decompose the underlying space into disjoint regions of attraction and establish the adjacency criterion. The criterion states a necessary and sufficient condition for two regions of attraction of stable equilibria to be adjacent. We then apply this criterion on a specific type of gradient flow which has as many as $n!$ local maxima. In particular, we characterize the set of equilibria, compute the index of each critical manifold and moreover, find all pairs of adjacent neighbors.  As an application of the adjacency criterion, we introduce a stochastic version of the double bracket flow and set up a Markov model to approximate the sample path behavior. The study of this specific prototype with its special structure provides insight into many other difficult problems involving simulated annealing.     \end{abstract}

\section{Introduction}

Ascent/Descent equations often provide the most direct demonstration of the existence of a maximum/minima and can provide an easily implemented algorithm to find the maximum/minima (see, for example,~\cite{brockett1991dynamical,bloch1990steepest,brockett1991gradient}).  Of course when the function being maximized/minimized has multiple local maxima/minima, steepest ascent/descent needs to be modified and for the last several decades the modification of choice has been some type of simulated annealing procedure.  However, because simulated annealing is slow and subject to variable results because of its stochastic nature, there remains considerable interest in finding methods for improving its speed and, in general, learning more about its performance. But often this requires the knowledge of the sample path behavior.

In this paper, we start with the development  of a basic theorem about the geometry of a general class of gradient flows. In particular,  we decompose the underlying space into disjoint regions of attraction associated with the gradient flow, and then establish the adjacency criterion. This theorem  states a necessary and sufficient condition for two regions of attraction of stable equilibria to be adjacent, i.e, for them to share a boundary of co-dimension one.

The adjacency criterion has a large potential impact on studying stochastic gradients.   For example, to approximate the sample path behavior of a stochastic gradient, R. W. Brockett established a Markov model in~\cite{brockett2011modeling} whose states consist of all stable equilibria. Each transition probability is evaluated by solving a related first-hitting time problem, the computational complexity is in a large-scale. The adjacency criterion  then sheds light on the problem.  Knowing adjacent neighbors reduces the computational amount as the transition probability between non-adjacent  equilibria is negligible when compared with the transition probability between adjacent equilibria, especially under the case where the noise is moderate.  

The adjacency criterion also relates to a number of geometric factors. For example, the depth of the potential well associated with a stable equilibrium, the volume of a region of attraction, the area of the boundary shared by a pair of adjacent neighbors and etc, all these factors are important for explaining simulated annealing. 
 


As an demonstration of adjacency criterion, we consider a prototype system with the gradient flow  of the type $\dot{H}=[H,[H,\pi(H)]]$, the map $\pi$ projects a matrix onto its diagonal.  This flow is a prototype for gradient systems with multiple stable equilibria, it has as many as $n!$ stable equilibria, each is a diagonal matrix and one-to-one corresponds to an element  in the permutation group $S_n$. In this paper, we will characterize all pairs of adjacent neighbors, the characterization is simple and clean. In particular, we will show that each pair of adjacent neighbors is  related to a simple transposition, so then the $n!$ regions of attraction are arranged in a way that each one has $(n-1)$ adjacent neighbors.  A stochastic version of the gradient flow is studied at the end of this paper as a concrete example of the application of the adjacency criterion. In particular, we set up an optimal control problem as an approach to evaluate the transition probability.

After this introduction, we proceed in steps: 
in section 2, we will introduce the adjacency criterion after a quick review of some basic notions about Morse-Bott gradient systems. In section 3, we will introduce the isospectral manifold, the normal metric and the double bracket flow under a special class of potentials functions. In section 4, we will focus on a specific quadratic potential function and show that the set of equilibria under the potential function is a mix of isolated points and continuum manifolds. In section 5, we will show that the quadratic potential is a Morse-Bott function by explicitly working out the Hessian of the potential. In particular, we will compute the index and the co-index of each critical manifold. In section 6, we will apply the adjacency criterion to characterize all pairs of adjacent neighbors. In the last section, we will work on a stochastic gradient as an application of the adjacency criterion. In particular, we will show how the adjacency criterion will simplify the evaluation of transition probability.

\section{Adjacency criterion}
In this section, we will introduce the adjacency criterion. Some mathematical backgrounds are needed here, in this paper we will only introduce terminologies that are necessary for establishing the criterion. 

A potential $\Psi$ on a Riemannian manifold $M$ is a \textbf{Morse-Bott function} if the set of equilibria under the gradient flow grad$(\Psi)$ is a finite disjoint union of connected submanifolds $\{E_1,\cdots,E_n\}$, and the Hessian of $\Psi$ is nondegenerate when restricted at the normal space $N_pE_k$ at any point $p\in E_k$ and for any $E_k$. (see, for example, \cite{banyaga2010morse} for more details about Morse-Bott functions) 

For convenience, we say each set $E_k$ is a \textbf{critical manifold}. Assume $\Psi$ is a Morse-Bott function, the \textbf{index} of each submanifold $E_k$ is the number of negative eigenvalues of the Hessian restricted at $N_p E_k$ for some $p\in E_k$, this is well-defined because the index of $E_k$ is independent of the choice of $p$. Similarly, the \textbf{co-index} of $E_k$ is the number of positive eigenvalues of the Hessian restricted at $N_p E_k$. An equation relating the index, co-index and the dimension of $M$ is  
\begin{equation}
\text{ind} E_k+\text{co-ind} E_k+\dim E_k=\dim M
\end{equation}
A {stable critical manifold} is then a critical manifold of co-index $0$.

Let $E_k$ be a critical manifold, the \textbf{stable manifold of $E_k$} is defined by
\begin{equation}
W^s(E_{k})=\{p\in M|\lim_{t\to\infty}\varphi_{t}(p)\in E_k \}
\end{equation}  
where $\varphi_t(p)$ is the solution of the differential equation $\dot p=$grad$(\Psi)$ parametrized by time $t$. The \textbf{unstable manifold of $E_k$} is defined in a similar way as
\begin{equation}
W^u(E_k)=\{p\in M|\lim_{t\to-\infty}\varphi_{t}(p)\in E_k\}
\end{equation}
The dimensions of $W^s(E_{k})$ and $W^u(E_{k})$ are ind$ E_k+\dim E_k$ and co-ind$E_k+\dim E_k$ respectively. A decomposition of $M$, with respect to stable/unstable manifolds of $E_k$'s, is given by
\begin{equation}
M=\bigcup^n_{k=1}W^s(E_k)
\end{equation}
or
\begin{equation}
M=\bigcup^n_{k=1}W^u(E_k)
\end{equation}
If only stable critical manifolds $\{S_1,\cdots,S_l\}$ are concerned, then 
\begin{equation}
M=\bigcup^l_{i=1}\overline{W^s(S_i)}
\end{equation}
We are interested in how these cells are arranged in the underlying space. With terminologies above, we are now ready to state the adjacency criterion. 
\vspace{5pt}

\begin{fact}[Adjacency Criterion] Suppose $\Psi$ is a Morse-Bott function on a smooth manifold $M$. Let $\{S_1,\cdots, S_l\}$ be the collection of stable critical manifolds and let $\{K_1,\cdots,K_m\}$ be the critical manifolds of co-index $1$. If 
\begin{equation}\label{ASSU}
\bigcup^m_{i=1}(W^u(K_i)-K_i)\subset \bigcup^l_{j=1}W^s(S_j)  
\end{equation}
then the boundary of each $W^s(S_j)$ is piece-wise smooth, each piece is a stable manifold of some $K_i$, i.e,
\begin{equation}
\partial\overline{W^s(S_j)}=\bigcup_{K_i\in \partial\overline{W^s(S_j)}}\overline{W^s(K_i)}
\end{equation}
So two stable manifolds $S_j$ and $S_k$ are adjacent if and only if there is a $K_i$ such that 
$
K_i \subset \partial\overline{W^s(S_j)}\cap\partial\overline{W^s(S_k)}
$.
\end{fact}
\vspace{5pt}

\begin{figure}[h]
\begin{center}
\includegraphics[scale=.4]{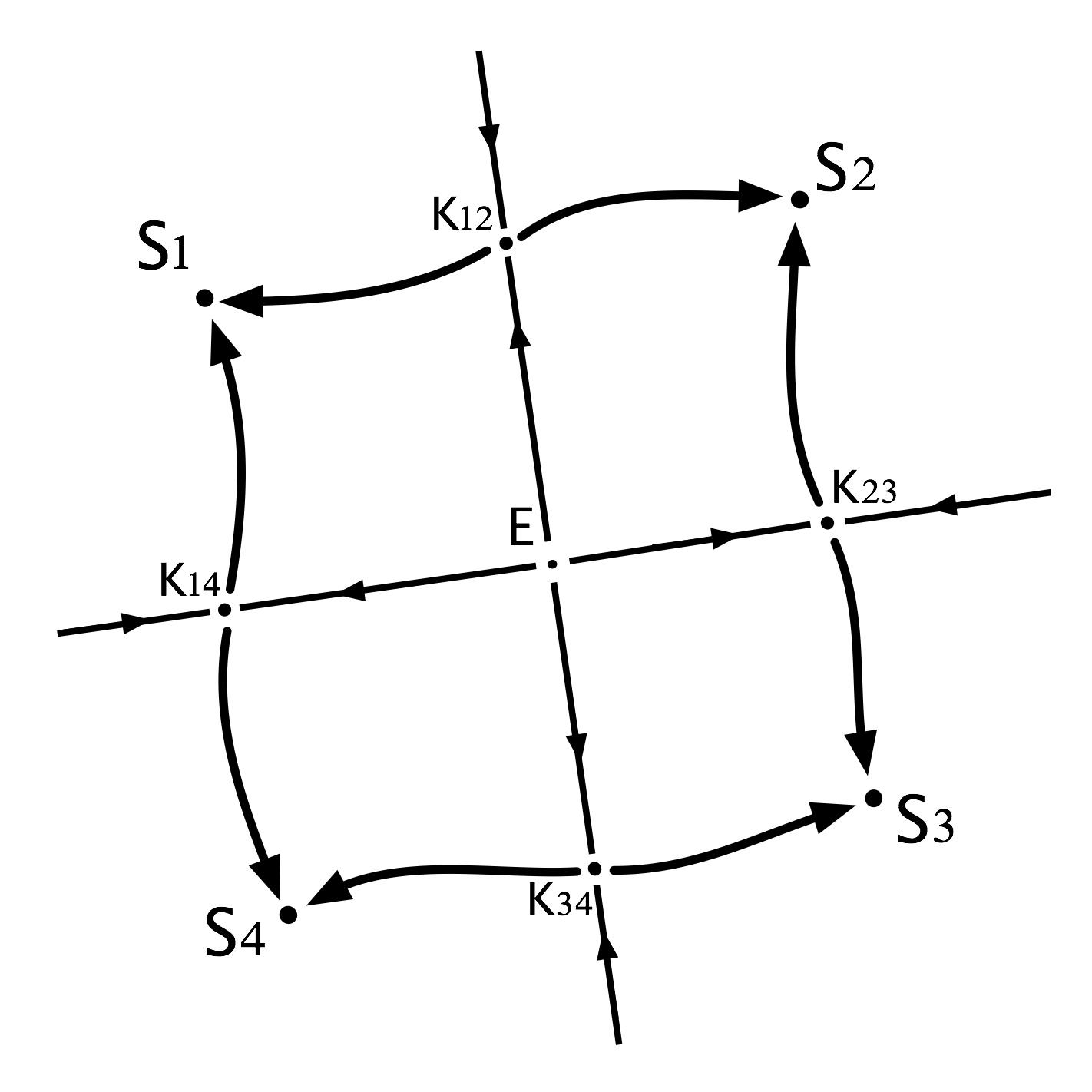}
\caption {In this figure, each $S_i$ is a stable equilibrium, each $K_{ij}$ is an equilibrium of co-index $1$ and $E$ is an equilibrium of co-index $2$. The pair  $(S_1, S_2)$, for example, consists of two adjacent neighbors because their regions of attraction  share the stable manifold of $K_{12}$ as a co-dimension one boundary. On the other hand, $S_1$ and $S_3$ are not adjacent because  their regions of attraction intersect each other only at $E$.}
\label{FIG1}
\end{center}
\end{figure}

Before heading to its application on the isospectral manifold, we point out that the assumption described by equation \eqref{ASSU} includes a broad range of gradient systems. For example, in a Morse-Smale gradient system, all equilibria are isolated points and the unstable manifold of each $K_i$ is a one dimensional curve and its boundaries consists of stable equilibrium points as illustrated in Fig. \ref{FIG1}. The class of Morse-Smale gradient system is a residual class in all gradient systems(see, for example, \cite{banyaga2013lectures,palis1969morse} for more details about Morse-Smale systems).

\section{Isospectral manifold and a class of symmetric potentials}
In this section, we will introduce the isospectral manifold with strongly disjoint eigenvalues, the normal metric and a special class of symmetric potential functions.

Let $Sym(\Lambda)$ denote the isospectral manifold, in case here, it refers to the collection of all $n$-by-$n$ real symmetric matrices with a fixed set of eigenvalues $\Lambda=\{\lambda_1,\cdots,\lambda_n\}$. We assume that the set $\Lambda$ is \textbf{strongly disjoint}, i.e, 
\begin{itemize}
\item for any two nonempty disjoint subsets $\{\lambda_{a_1},\cdots,\lambda_{a_p}\}$ and $\{\lambda_{b_1},\cdots,\lambda_{b_q}\}$ of $\Lambda$, the inequality $\frac{1}{p}\sum^p_{k=1}\lambda_{a_k}\neq\frac{1}{q}\sum^q_{k=1}\lambda_{b_k}$ holds. 
\end{itemize}
The assumption on $\Lambda$ is a stronger version of the condition that all eigenvalues are distinct. Yet this condition is generic in the sense that if we randomly pick a vector $(\lambda_1,\cdots,\lambda_n)$ in $\mathbb{R}^n$, then $\Lambda$ is strongly disjoint for almost sure because the set defined by equations $\frac{1}{p}\sum^p_{k=1}\lambda_{a_k}=\frac{1}{q}\sum^q_{k=1}\lambda_{b_k}$ is a finite union of hyperplanes in $\mathbb{R}^n$. On the other side, as we will see later this condition on $\Lambda$ is necessary and sufficient for a class of potential functions on $Sym(\Lambda)$ to be Morse-Bott functions which is a key assumption for application of the adjacency criterion. 

The so called \textbf{normal metric} $g$ on $Sym(\Lambda)$ is defined as follows, the tangent space $T_HSym(\Lambda)$ consists of elements of the form $[H,\Omega]$ with $\Omega$ skew symmetric, since all eigenvalues are distinct, the adjoint map $ad_H:\Omega\mapsto [H,\Omega]$ is then an isomorphism. For any two tangent vectors $[H,\Omega_1]$, $[H,\Omega_2]$ in $T_HSym(\Lambda)$, the normal metric $g$ at $H$ is then defined by
\begin{equation}
g([H,\Omega_1],[H,\Omega_2]):=-tr(\Omega_1\Omega_2)
\end{equation}
It is routine to check that $g$ is positive definite. Equipped with the normal metric, the gradient flow of a smooth function $\Psi\in C^{\infty}(Sym(\Lambda))$ is then the double bracket flow
\begin{equation}
\dot{H}=[H,[H,\Psi'(H)]]
\end{equation}
where $\Psi'(H)$ denotes the derivative of $\Psi$ with respect to $H$

For convenience, we assume that eigenvalues are ordered as $\lambda_1<\cdots<\lambda_n$ and we denote by $d_1,\cdots,d_n$ the diagonal entries of a matrix $H$ in $Sym(\Lambda)$. There is a special class of potentials functions $\Psi$ on $Sym(\Lambda)$, each is symmetric in diagonal entries and can be generated by a scalar function $\phi\in C^{\infty}[\lambda_1,\lambda_n]$ via the equation 
\begin{equation}\label{POT}
\Psi(d_1,\cdots,d_n)=\sum^n_{i=1}\int^{d_i}_{\lambda_1}\phi(x)dx
\end{equation}   
As the sum of the diagonal entries is constant for any matrix in $Sym(\Lambda)$, there is actually a broad class of symmetric potentials can be expressed in this way. In fact, we have shown in~\cite{chen2014symmetric} that if $\Psi(d_1,\cdots,d_n)$ can be expressed as a power series of symmetric polynomials, then there is a scalar function $\phi$ such that the equation \eqref{POT} holds. Moreover, we have also shown in \cite{chen2014symmetric} that for almost all scalar functions $\phi$, the potential $\Psi$ generated by equation \eqref{POT} is a Morse-Bott function.

In this paper, we will work on a simple case where $\phi(x)=x$. The corresponding potential is then in a diagonal form as
\begin{equation}\label{POT1}
\Psi=\frac{1}{2}\sum^n_{i=1}d^2_i
\end{equation}
In the rest of paper, we will just call this specific quadratic potential the \textbf{diagonal potential}. In the next two sections, we will characterize the critical manifolds associated with the gradient flow and compute explicitly the Hessian at each equilibrium point.

\section{The set of equilibria associated with the diagonal potential}
For convenience, we let $\pi$ be a projection map sending each matrix $H$ to the diagonal matrix $diag(d_1,\cdots,d_n)$. Then the gradient vector field $f(H)$ with respect to the diagonal potential is simply 
\begin{equation}
f(H)=[H,[H,\pi(H)]]
\end{equation}
A matrix $H\in Sym(\Lambda)$ is an equilibrium if $f(H)=0$ and it happens if and only if $[H,\pi(H)]=0$ because $tr(\pi(H)f(H))=-tr([H,\pi(H)]^2)$. This leads us to   
\vspace{5pt}

\begin{lem}\label{LBLO} If $H$ is an equilibrium, then there exists a permutation matrix $P$ such that $P^THP$ is block-diagonal, i.e,
\begin{equation}\label{BLOH}
P^THP=\text{Diag}(H_1,\cdots,H_k)
\end{equation} 
Suppose $d_{i1},\cdots,d_{in_i}$ are diagonal entries of $H_i$, then  
\begin{equation}\label{DEFPHII}
d_{i1}=\cdots=d_{in_i}
\end{equation}
and this holds for each $H_i$. 
\end{lem}
\vspace{5pt}

\begin{proof} The commutator $[H,\pi(H)]$ vanishes if and only if $h_{ij}(d_i-d_j)=0$ for each pair of $(i,j)$.   
\end{proof} 
\vspace{5pt}

Let $\Lambda_i$ be the set of eigenvalues of $H_i$, then $(\Lambda_1,\cdots,\Lambda_k)$ is \textbf{a partition of $\Lambda$}, i.e, $\Lambda=\bigcup^k_{i=1}\Lambda_i$ and $\Lambda_i\cap\Lambda_j=\varnothing$  if $i\neq j$. For convenience, we let $\#(\Lambda_i)$ denote the cardinality of $\Lambda_i$, let $s(\Lambda_i)$ denote the sum of its elements and we define
\begin{equation}
\mu(\Lambda_i):=s(\Lambda_i)/\#(\Lambda_i)
\end{equation}
then 
\begin{equation}
d_{i1}=\cdots=d_{in_i}=\mu(\Lambda_i)
\end{equation}

A symmetric matrix $H$ is said to be \textbf{irreducible} if there isn't a permutation matrix $P$ such that $P^TH P$ is a nontrivial block-diagonal matrix.   
\vspace{5pt}
 
\begin{lem} Each $H_i$ in equation \eqref{BLOH} is irreducible, moreover $\mu(\Lambda_i)\neq \mu(\Lambda_j)$ if $i\neq j$.  
\end{lem}
\vspace{5pt}

\begin{proof} Suppose $H_i$ is not irreducible, we may assume that $H_i=Diag(H_{i1},H_{i2})$, let $\Lambda_{i1}=\Lambda_{i2}$ be the set of eigenvalues of $H_{i1}$ and $H_{i2}$ respectively, then $\mu(\Lambda_{i1})=\mu(\Lambda_{i2})=\mu(\Lambda_i)$ which contradict the fact that $\Lambda$ is strongly disjoint. By applying the same arguments, we conclude that $\mu(\Lambda_i)$ and $\mu(\Lambda_j)$ are disjoint if $i\neq j$. 
\end{proof}
\vspace{5pt}

An equilibrium $H$ gives rise to a partition of $\Lambda$. Conversely a partition of $\Lambda$ will correspond to a set of equilibria. Let $\mathcal{A}$ be the collection of all choices of partitions of $\Lambda$. Given a choice of partition $\alpha=(\Lambda_1,\cdots,\Lambda_k)$, we define a subset of $Sym(\Lambda)$ as
\begin{equation}
E_{\alpha}:=\{Diag(H_1,\cdots,H_k)|\pi(H_i)=\mu(\Lambda_i)\}
\end{equation}
clearly $E_{\alpha}$ consists exclusively of equilibria. Moreover it is a smooth manifold as a consequence of the next fact
\vspace{5pt}

\begin{fact}\label{Fact1} Let $\Lambda'=\{\lambda'_1,\cdots,\lambda'_{n'}\}$ be a subset of $\Lambda$ and let $C(\Lambda')\subset\mathbb{R}^{n'}$ be the convex hull of all vectors $(\lambda'_{\sigma(1)},\cdots,\lambda'_{\sigma(n')})$ where $\sigma$ varies over all permutations of $\{1,\cdots,n'\}$, then the image of the projection
\begin{equation}
\pi:Sym(\Lambda')\longrightarrow \mathbb{R}^{\#(\Lambda')}
\end{equation} 
is the convex hull $C(\Lambda')$(see, for example, \cite{horn1954doubly}). Let $\vec e$ be a vector of all ones in $\mathbb{R}^{\#(\Lambda')}$ and define
\begin{equation}
X(\Lambda'):=\pi^{-1}(\mu(\Lambda')\vec e)
\end{equation}
then $\pi$ is a submersion at each point $H\in X(\Lambda')$.
\end{fact}  
\vspace{5pt}

\begin{figure}[h]
\begin{center}
\includegraphics[scale=.4]{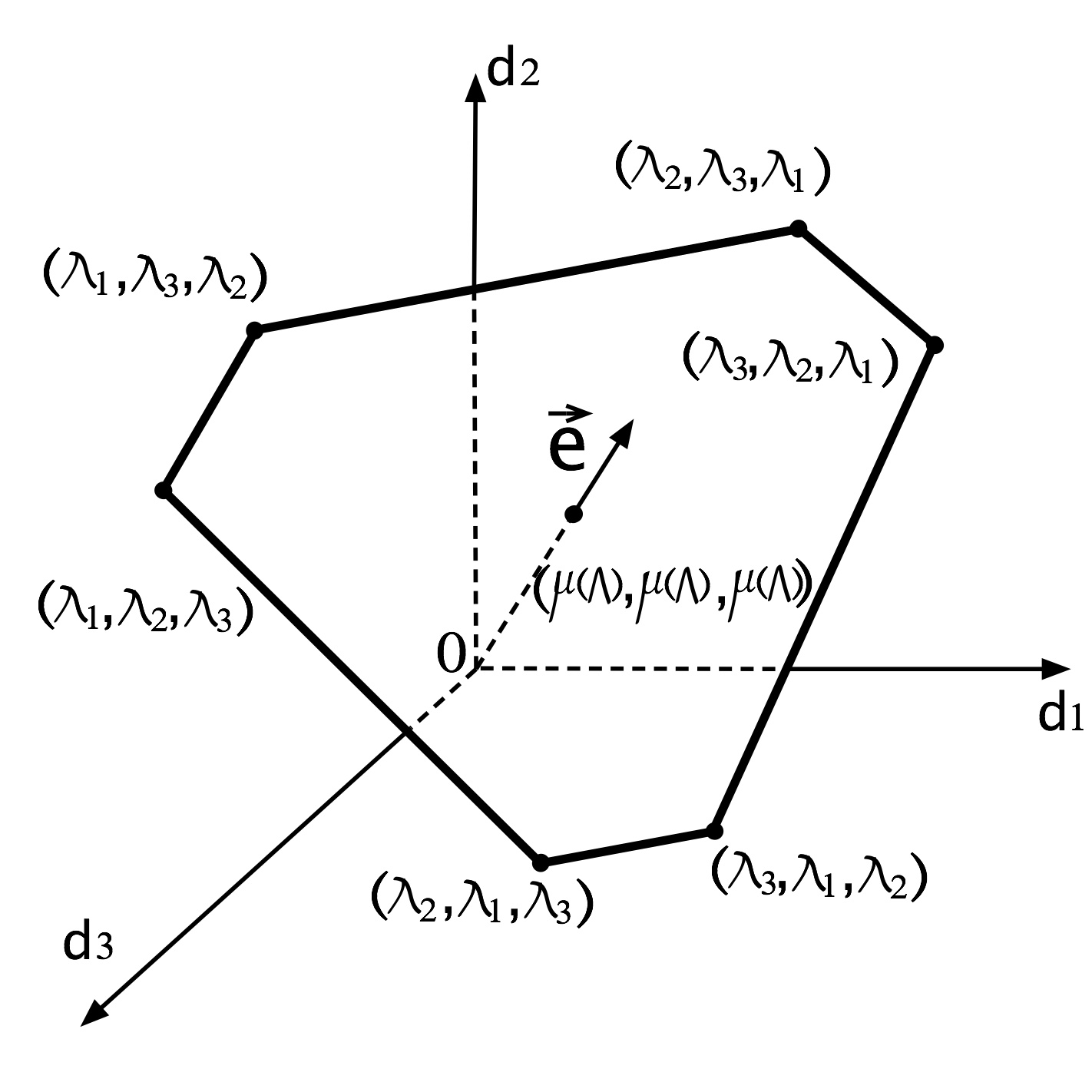}
\caption {The convex hull $C(\Lambda)$ under the case $\Lambda=\{\lambda_1,\lambda_2,\lambda_3\}$. Let $\vec e:=(1,1,1)^T$, then the line $\mathbb{R}\vec e$ intersects $C(\Lambda)$ perpendicularly at the point $\mu(\Lambda)\vec e.$}
\end{center}
\end{figure}

A complete proof can be found in \cite{chen2014symmetric}. Assuming the fact above, then an immediate consequence is that $X(\Lambda')$ is a smooth manifold of dimension $\frac{1}{2}(\#(\Lambda')-1)(\#(\Lambda')-2)$. An identification is that  $E_{\alpha}\simeq\Pi^k_{i=1}X(\Lambda_i)$, so 
\begin{equation}
\dim E_{\alpha}=\sum^k_{i=1}(\#(\Lambda_i)-1)(\#(\Lambda_i)-2)
\end{equation}  
An upshot is that the set $E_{\alpha}$ is a set of discrete points if and only if each $\#(\Lambda_i)$ is either one or two, in the case all $\Lambda_i$'s are singletons, the set $E_{\alpha}$ consists only of a diagonal matrix.   

At this moment, we have characterized the set equilibria that are block-diagonal, equilibria that are off block-diagonal can be generated by letting the group of permutation matrices act on $E_{\alpha}$'s by conjugation as we will describe below.      
\vspace{5pt}

\begin{lem}\label{LEMSEC1} If $H$ is an equilibrium, then so is $P^THP$ for any permutation matrix $P$. 
\end{lem}
\vspace{5pt}

\begin{proof} We check that 
$
\pi(P^THP)=P^T\pi(H)P
$
and hence conjugation by permutation matrices commutes with the commutator, i.e, 
\begin{equation}
P^T[H,\pi(H)]P=[P^THP, \pi(P^THP)]
\end{equation}
So $f(P^THP)=0$ exactly when $f(H)=0$.
\end{proof}
\vspace{5pt}

Lemma \ref{LEMSEC1} can be regarded as a converse argument of lemma \ref{LBLO}. 
\vspace{5pt}

\begin{theorem} Suppose $\alpha=(\Lambda_1,\cdots,\Lambda_k)$ is a choice of partition and suppose $\#(\Lambda_i)=n_i$. Let the group of permutation matrices act on $E_{\alpha}$ by conjugation, then the orbit of $E_{\alpha}$ contains as many as $n!/\Pi^k_{i=1}n_i!$ disjoint smooth submanifolds in $Sym(\Lambda)$, each consists exclusively of equilibria. The set of equilibria is then the union of orbits of $E_{\alpha}$ as $\alpha$ varies over $\mathcal{A}$. 
\end{theorem}
\vspace{5pt}

\begin{proof} For any permutation matrix $P$, the set $P^TE_{\alpha}P$ consists exclusively of equilibria by lemma \ref{LEMSEC1}. On the other hand, either  
$
P^TE_{\alpha}P=E_{\alpha}
$
or 
$
P^TE_{\alpha}P\cap E_{\alpha}=\varnothing
$, 
and $P$ fixes $E_{\alpha}$ if and only if $P$ is block-diagonal, i.e, 
$
P=Diag(P_1,\cdots,P_k)
$
with each $P_i$ a $n_i$-by-$n_i$ permutation matrix. There are exactly $\Pi^k_{i=1}n_i!$ many block-diagonal permutations matrices, so by Burnside's counting formula, there are as many as $n!/\Pi^k_{i=1}n_i!$ disjoint sets in the orbit of $E_{\alpha}$, each is a smooth submanifold in $Sym(\Lambda)$. 
\end{proof}
\vspace{5pt}

In the rest of this paper, we will abuse the term \textit{critical manifold} by referring it to an $E_{\alpha}$ or any set in its orbit as $E_{\alpha}$ may have multiple connected components. (see, for example, a discussion on the number of components in \cite{chen2014symmetric}).  Notice that in the orbit of $E_{\alpha}$, there exist multiple critical manifolds that are block-diagonal, for example, the orbit of a diagonal matrix is the set of all diagonal matrices. Actually given a choice of partition $\alpha=(\Lambda_1,\cdots,\Lambda_n)$, there are as many as $k!$ critical manifolds that are block-diagonal leaving the others off block-diagonal. If we want to choose a \textbf{canonical representative} among the orbit, then permute $\Lambda_1,\cdots,\Lambda_k$ if necessary so that $\mu(\Lambda_1)<\cdots<\mu(\Lambda_k)$, this is possible because eigenvalues of $\Lambda$ are strongly disjoint.

\section{The Hessian of the diagonal potential} 
In this section, we will show that the diagonal potential is a Morse-Bott function by explicitly working out the Hessian of the potential at each equilibrium point and compute its eigenvalues.  

Let $\mathcal{H}$ denote the Hessian of a potential $\Psi$. On a Riemannian manifold, it is defined by 
\begin{equation}
\mathcal{H}=\nabla^2\Psi
\end{equation} 
where $\nabla$ is the Levi-Civita connection and if we evaluate the Hessian at an equilibrium point, then 
\begin{equation}\label{HESSFOR}
\mathcal{H}(X,Y)=X(g(\text{grad}(\Psi),Y))
\end{equation}
In this specific case, we have 
\vspace{5pt}

\begin{fact}
Suppose $H$ is an equilibrium, then the Hessian $\mathcal{H}$ evaluated at $H$ is given by
\begin{eqnarray}\label{HES}
\mathcal{H}([H,\Omega_i],[H,\Omega_j]) = & -tr([H,\Omega_i][\pi(H),\Omega_j])\notag\\
 &-\langle\pi([H,\Omega_i]),\pi([H,\Omega_j])\rangle
\end{eqnarray}
where $\langle\cdot,\cdot\rangle$ is the normal inner-product in $\mathbb{R}^n$. 
\end{fact}  
\vspace{5pt}

We omit the proof here as this is a direct computation following equation \eqref{HESSFOR}. Given a choice of partition $\alpha=(\Lambda_1,\cdots,\Lambda_n)$, we assume $\#(\Lambda_i)=n_i$, then the dimension of the normal space at a point $H\in E_{\alpha}$ is given by
\begin{equation}
\dim N_H E_{\alpha}=\frac{1}{2}\{n(n-1)-\sum^k_{i=1}(n_i-1)(n_i-2)\}
\end{equation}
For convenience, let $n_{\alpha}:=\dim N_H E_{\alpha}$. We will now construct an orthogonal basis $\mathcal{N}$ of $N_H E_{\alpha}$ with respect to $\mathcal{H}$, i.e, 
\begin{eqnarray}\label{HECR}
& \mathcal{H}(N,N)\neq 0, \forall N\in\mathcal{N}\\
& \mathcal{H}(N,N')=0, \text{ if } N\neq N'
\end{eqnarray} 

Suggested by the partition $\alpha$, we divide a $n$-by-$n$ matrix $N$ into $k$-by-$k$ blocks as 
\begin{equation}\label{DIVI}
N=
\begin{pmatrix}
B_{11} & \cdots & B_{1k}\\
\vdots & \ddots & \vdots\\
B_{k1} & \cdots & B_{kk}
\end{pmatrix}
\end{equation}
with the $pq$-th block of dimension $n_p$-by-$n_q$. The basis $\mathcal{N}$ consists of two parts: the block-diagonal part $\mathcal{N}_d$ and the off block-diagonal part $\mathcal{N}_o$.
\vspace{5pt}

\subsubsection*{1. Constructing $\mathcal{N}_o$}
Write $H=Diag(H_1,\cdots,H_k)$ with each $H_i\in Sym(\Lambda_i)$ and let $\vec v_{i1},\cdots,\vec v_{in_i}$ be the unit-length eigenvectors of $H_i$ with respect to the eigenvalues $\lambda_{i1},\cdots,\lambda_{in_i}$.

Given integers $p$ and $q$, define a set of symmetric matrices $\mathcal{N}_{pq}$ in the way that for each matrix $N_{pq,ij}\in \mathcal{N}_{pq}$, all blocks except $B_{pq}$ and $B_{qp}$ are zeros while
\begin{equation}
B_{pq}(=B^T_{qp})=(\lambda_{pi}-\lambda_{qj})\vec v_{pi}\vec v^T_{qj}
\end{equation}
for some $i=1,\cdots,n_p$ and some $j=1,\cdots,n_q$, so there are exactly $n_pn_q$ symmetric matrices in the set $\mathcal{N}_{pq}$. 
\vspace{5pt}

\begin{lem}\label{OFFDIAG} The set 
$
\mathcal{N}_o:=\bigcup_{p<q}\mathcal{N}_{pq}
$ 
is contained in $N_{H}E_{\alpha}$. If $N_{pq,ij}\in \mathcal{N}_{pq}$ and $N_{p'q',i'j'}\in \mathcal{N}_{p'q'}$, then
\begin{align}\label{INDE}
& \mathcal{H}(N_{pq,ij},N_{p'q',i'j'})\notag\\
= & -2(\lambda_{pi}-\lambda_{qj})(\mu(\Lambda_p)-\mu(\Lambda_q))\Delta_{pp',qq',ii',jj'}
\end{align}
with $\Delta_{pp',qq',ii',jj'}:=\delta_{pp'}\delta_{qq'}\delta_{ii'}\delta_{jj'}$ and $\delta_{ij}$ is the Kronecker delta.  
\end{lem}    
\vspace{5pt}

This is a straitforward computation following the formula \eqref{HES}, more details can be found in \cite{chen2014symmetric}. The set $\mathcal{N}_o$ contains as many as $\sum_{1\le i<j\le k}n_in_j$ off-block-diagonal matrices. 
\vspace{5pt}

\subsubsection*{2. Constructing $\mathcal{N}_d$} still assume $H=Diag(H_1,\cdots,H_k)$ with each $H_i\in Sym(\Lambda_i)$. Fix $i$, we let $\vec e\in \mathbb{R}^{n_i}$ be a vector with all entries ones and let $\vec e^{\perp}$ be the hyperplane perpendicular to $\vec e$. An identification is that the tangent space of the convex hull $C(\Lambda')$ at any of its point is $\vec e^{\perp}$. 
\vspace{5pt}

\begin{lem} Let $\{\vec u_1,\cdots, \vec u_{n_i-1}\}$ be an orthonormal basis of $\vec e^{\perp}$, then there exists a set of skew symmetric matrices $\{\tilde \Omega_{i,1},\cdots,\tilde\Omega_{i,n_i-1}\}$ such that  
$
\pi([H_i,\tilde\Omega_{i,j}])=u_j
$.
\end{lem}
\vspace{5pt}

\begin{proof} This follows fact \ref{Fact1} as the projection $\pi$ is a submersion. 
\end{proof}
\vspace{5pt}

Define a set of symmetric matrices $\mathcal{N}_i$ in the way that for each matrix $N_{i,j}\in\mathcal{N}_i$, all blocks except $B_{ii}$ are zeros while 
\begin{equation}
B_{ii}=[H_i,\tilde\Omega_{i,j}] 
\end{equation}
for some $j=1,\cdots,n_i-1$, so there are exactly $(n_i-1)$ symmetric matrices in the set $\mathcal{N}_i$.   
\vspace{5pt}

\begin{lem}\label{DIAGO} The set 
$
\mathcal{N}_d:=\bigcup^k_{i=1}\mathcal{N}_i
$ 
is contained in $N_HE_{\alpha}$. If $N_{i,j}\in\mathcal{N}_{i}$ and $N_{i',j'}\in\mathcal{N}_{i'}$, then
$
\mathcal{H}(N_{i,j},N_{i',j'})=\delta_{ii'}\delta_{jj'}
$.  
\end{lem}
\vspace{5pt}

This again follows the formula \eqref{HES}. The set $\mathcal{N}_d$ contains as many as $\sum^k_{i=1}(n_i-1)$ block-diagonal matrices. 
\vspace{5pt}

\subsubsection*{3. Orthogonality of $\mathcal{N}_o$ and $\mathcal{N}_d$} 

First notice that by combining $\mathcal{N}_o$ and $\mathcal{N}_d$, there are exactly $n_{\alpha}$ matrices, i.e,  
\begin{equation}
n_{\alpha}=\sum_{1\le i<j\le k}n_in_j+\sum^k_{i=1}(n_i-1)
\end{equation} 
because of the equality $\sum^k_{i=1}n_i=n$. So it all remains to show that 
\vspace{5pt}

\begin{lem}\label{ORTHO} Matrices in $\mathcal{N}_o$ are orthogonal to those in $\mathcal{N}_d$ with respect to the Hessian. 
\end{lem}
\vspace{5pt}

\begin{proof} Divide any $n$-by-$n$ skew symmetric matrix $\Omega$ into $k$-by-$k$ blocks as we did in equation \eqref{DIVI}. Define a skew symmetric matrix $\Omega_{pq,ij}$ by setting all blocks but $B_{pq}$ and $B_{qp}$ zeros while
\begin{equation}
B_{pq}(=-B^T_{qp})=\frac{1}{2}\vec v_{pi}\vec v^T_{qj} 
\end{equation} 
then it is a straitforward computation that 
\begin{equation}
N_{pq,ij}=[H,\Omega_{pq,ij}]
\end{equation}
On the other side, if we define a skew symmetric matrix $\Omega_{s,t}$ by setting all blocks zeros but leaving $B_{ss}=\tilde\Omega_{s,t}$, then  
\begin{equation}
N_{s,t}=[H,\Omega_{s,t}]
\end{equation} 
A computation shows that 
\begin{eqnarray}
& [\pi(H), \Omega_{s,t}]=0\\
& \pi([H,\Omega_{pq,ij}])=0
\end{eqnarray}
These two equalities annihilate the right hand side of formula \eqref{HES}, so $\mathcal{H}(N_{pq,ij},N_{s,t})=0$.   
\end{proof}
\vspace{5pt}
 
\begin{theorem} The Hessian $\mathcal{H}$ is nondegenerate when restricted at the normal space $N_H E_{\alpha}$ for any $H\in E_{\alpha}$ and any $E_{\alpha}$. Moreover, $\mathcal{H}$ is invariant under conjugation of permutation matrices, so the diagonal potential by equation $\eqref{POT1}$ is a Morse-Bott function. 
\end{theorem}
\vspace{5pt}

\begin{proof} The invariance is an consequence of lemma \ref{LEMSEC1} and formula \eqref{HES}. The rest is an outcome by combining lemma \ref{OFFDIAG}, lemma \ref{DIAGO} and lemma \ref{ORTHO}.
\end{proof}
\vspace{5pt}

We end this section with a discussion on the basis $\mathcal{N}$. For each symmetric matrix $N\in N_H E_{\alpha}$, there is a unique matrix $N'\in N_H E_{\alpha}$ such that $\mathcal{H}(N,\cdot)=tr(N'\cdot)$ because $\mathcal{H}$ is nondegenerate. This induces a linear map $L_{\mathcal{H}}$ on $N_H E_{\alpha}$ by sending $N$ to $N'$. Each linear subspace spanned by $\mathcal{N}_{pq}$ or $\mathcal{N}_{i}$ is invariant under $L_{\mathcal{H}}$, moreover each matrix $N_{pq,ij}\in \mathcal{N}_{pq}$ is an eigenmatrix of $L_{\mathcal{H}}$ with respect to the eigenvalue $-(\mu(\Lambda_p)-\mu(\Lambda_q))/(\lambda_{pi}-\lambda_{qj})$. All eigenvalues from the linear subspace spanned by $\mathcal{N}_d$ are positive by lemma \ref{DIAGO}. So at this moment for each critical manifold, we have computed its index and co-index.

\section{Application of adjacency criterion on the isospectral manifold with the diagonal potential} 
In this section, we will work out all stable critical manifolds and characterize all pairs of adjacent neighbors. 
\vspace{1pt}

\begin{lem}\label{STABEQ} A critical manifold is stable if and only if it is a singleton consisting of a diagonal matrix.
\end{lem}
\vspace{5pt}

\begin{proof}
Suppose $M$ is a stable critical manifold, without loss of generality, we assume $M=E_{\alpha}$ because the Hessian is invariant under conjugation of permutation matrices, if $E_{\alpha}$ is stable, then so is any critical manifold in its orbit. Suppose $\alpha=(\Lambda_1,\cdots,\Lambda_k)$, then by lemma \ref{DIAGO}, each $\Lambda_i$ is a singleton, otherwise there exist at least $\sum^k_{i=1}(\#(\Lambda_i)-1)$ positive eigenvalues of $L_{\mathcal{H}}$. So $M$ is necessarily a singleton consisting of a diagonal matrix. On the other hand, if $H$ is a diagonal matrix, then it is a stable equilibrium because all eigenvalues of $L_{\mathcal{H}}$ are equal $-(\lambda_{pi}-\lambda_{qj})/(\lambda_{pi}-\lambda_{qj})$=-1 as suggested by our discussion in the end of last section.  
\end{proof}
\vspace{5pt}

A diagonal matrix as an stable equilibrium is worth having its own notation $s_{\sigma}$, the subindex $\sigma$ indicates the permutation on the set of indices $\{1,\cdots,n\}$, i.e, $s_{\sigma}=(\lambda_{\sigma(1)},\cdots,\lambda_{\sigma(n)})$. We now characterize critical manifolds of co-index $1$. Recall that $\lambda_1<\cdots<\lambda_n$, we say that two eigenvalues $\lambda_i<\lambda_j$ are \textbf{close in order} if $j=i+1$. 
\vspace{5pt}

\begin{lem} Suppose $\alpha=(\Lambda_1,\cdots,\Lambda_k)$ and $E_{\alpha}$ is a critical manifold is of co-index $1$, then all but one $\Lambda_p$ are singletons, and $\Lambda_p$ consists of two eigenvalues and they are closed in order.    
\end{lem}
\vspace{5pt}

\begin{proof} By lemma \ref{DIAGO} and lemma \ref{STABEQ}, we first conclude that there is exactly one  $\Lambda_p$ among $\{\Lambda_1,\cdots,\Lambda_k\}$ such that it is not a singleton and $\#(\Lambda_p)=2$. The matrix set $\mathcal{N}_p$ is then a singleton contains exactly one block-diagonal matrix $N_p$. 

At this moment, there is at least one positive eigenvalue of $L_{H}$ with eigenmatrix $N_p$ because of lemma \ref{DIAGO} and the fact that span$\mathcal{N}_p$ is an invariant subspace. Suppose $\Lambda_p=\{\lambda_{p1},\lambda_{p2}\}$, then to prevent from having extra positive eigenvalues of $L_{\mathcal{H}}$, two eigenvalues $\lambda_{p1}$ and $\lambda_{p2}$ are necessarily closed in order because if not, suppose $\lambda_q$ is in between $\lambda_{p1}$ and $\lambda_{p2}$. We consider the two eigenvalues of $L_{\mathcal{H}}$ contributed by the two dimensional linear subspace spanned by $\mathcal{N}_{pq}$, these two eigenvalues are $-(\mu(\Lambda_p)-\lambda_q))/(\lambda_{p1}-\lambda_{q})$ and $-(\mu(\Lambda_p)-\lambda_q)/(\lambda_{p2}-\lambda_{q})$. By assumption one is negative while the other is positive. 

On the other hand, suppose $\lambda_{p1}<\lambda_{p2}$ and they are close in order, since $\lambda_{p1}<\mu(\Lambda_p)<\lambda_{p2}$, it is clear that $-(\mu(\Lambda_p)-\lambda_q))/(\lambda_{pi}-\lambda_{q})<0$ for any $i=1,2$ and any $\lambda_q\in\Lambda-\{\lambda_{i1},\lambda_{i2}\}$.                 
\end{proof}
\vspace{5pt}

Each critical manifold of co-index $1$ is a discrete set of two elements. A typical set looks like 
\begin{equation}\label{COIND1}
\frac{1}{2}
\begin{pmatrix}
\ddots & & & &\\
& \lambda_{i}+\lambda_{i+1}& \cdots & \pm (\lambda_{i}-\lambda_{i+1}) & \\
& \vdots& \ddots&\vdots &\\
&\pm (\lambda_{i}-\lambda_{i+1})& \cdots& \lambda_{i}+\lambda_{i+1}& \\
& & & & \ddots
\end{pmatrix}
\end{equation}
To relate with diagonal matrices, we let $S_n$ be the group of permutations on indices $\{1,\cdots,n\}$, in convention a permutation is said to be a \textbf{simple transposition} if it is a 2-cycle $(i,i+1)$. We say two permutations $\sigma_1$ and $\sigma_2$ are related by a simple transposition $\hat\sigma$ if $\sigma_1\cdot\hat\sigma=\sigma_2$. Suppose 
\begin{eqnarray}
s_{\sigma_1}=diag(\cdots,\lambda_{i},\cdots,\lambda_{i+1},\cdots)   \\
s_{\sigma_2}=diag(\cdots,\lambda_{i+1},\cdots,\lambda_{i},\cdots)   
\end{eqnarray}
we then denote by $K_{\sigma_1,\sigma_2}$ the two matrices in equation \eqref{COIND1}. The collection of equilibria of co-index $1$ is then the union of $K_{\sigma_1,\sigma_2}$ where $(\sigma_1,\sigma_2)$ varies over all pairs of permutations that are related by a simple transposition.      

A closed set $Z\subset Sym(\Lambda)$ is \textbf{invariant} if for any $H\in Z$, the solution $\varphi_t(H)$ remains in $Z$ for any $t\in\mathbb{R}$. An observation is that
\vspace{5pt}

\begin{lem}\label{INVA} 
Let  $\alpha=(\Lambda_1,\cdots,\Lambda_k)$ be a choice of partition. We define 
\begin{equation}
Z_{\alpha}:=\{Diag(H_1,\cdots,H_k)|H_i\in Sym(\Lambda_i)\}
\end{equation}
then $Z_{\alpha}$ is an invariant subset in $Sym(\Lambda)$. 
\end{lem}
\vspace{5pt}

\begin{proof} If $H\in Z_{\alpha}$ and we write $H=Diag(H_1,\cdots,H_k)$, then the dynamical system is decoupled in the sense that 
\begin{equation}
\dot H_i=[H_i,[H_i,\pi (H_i)]], \forall i=1,\cdots,k 
\end{equation}
and $\dot H=Diag(\dot H_1,\cdots,\dot H_k)$.
\end{proof}
\vspace{5pt}

\textit{Remark:} Since the gradient flow $f(H)$ commutes conjugation by permutation matrices, so $P^TZ_{\alpha}P$ is also an invariant subset of $Sym(\Lambda)$ for any permutation matrix $P$.  
\vspace{5pt}

\begin{theorem}[Adjacent neighbors of diagonal matrices] All stable critical manifolds are singletons, and they are diagonal matrices. Two diagonal matrices $s_{\sigma_1}$ and $s_{\sigma_2}$ are adjacent if and only if $\sigma_1$ and $\sigma_2$ are related by a simple transposition. 
\end{theorem}
\vspace{5pt}

\begin{proof} Fix the pair $(\sigma_1,\sigma_2)$, and assume $k_{+}$ and $k_{-}$ are the two matrices in $K_{\sigma_1,\sigma_2}$. Let $W^{u}(k_{+})$ and $W^u(k_{-})$ be the unstable manifolds of $k_+$ and $k_-$ respectively, we will show that the boundary of either of unstable manifolds is $\{s_{\sigma_1},s_{\sigma_2}\}$, i.e, 
\begin{eqnarray}
\partial \overline{W^u(k_+)} = \partial \overline{W^u(k_-)} = s_{\sigma_1}\cup s_{\sigma_2}
\end{eqnarray}
But this is true by lemma \ref{INVA} because if let $\Lambda':=\{\lambda_1,\lambda_2\}$ and consider the double bracket flow on $Sym(\Lambda')$, then there are exactly four isolated equilibria, two diagonal matrices as stable equilibria and the other two are 
\begin{equation}
k'_{\pm}=\frac{1}{2}
\begin{pmatrix}
\lambda_{i}+\lambda_{i+1} & \pm(\lambda_i-\lambda_{i+1}) \\
\pm(\lambda_i-\lambda_{i+1}) & \lambda_{i}+\lambda_{i+1} 
\end{pmatrix}
\end{equation}
The isospectral manifold $Sym(\Lambda')$ is diffeomorphic to the circle $S^{1}$, if we parametrize it by $\theta\in\mathbb{R}$,  then the induced dynamical system of $\theta$ on the covering space $\mathbb{R}$ is then
\begin{equation}\label{EXMP}
\dot\theta=-\frac{1}{2}(\lambda_{i}-\lambda_{i+1})^2\sin(2\theta)
\end{equation}
It is clear that $\mathbb{Z}\pi$ are stable equilibria and $(\mathbb{Z}+\frac{1}{2})\pi$ are the unstable ones, this is consistent with the earlier arguments. Moreover this implies that the boundary of either $\overline{W^u(k'_+)}$ or $\overline{W^u(k'_-)}$ is the union of $diag(\lambda_i,\lambda_{i+1})$ and $diag(\lambda_{i+1},\lambda_i)$. This then completes the proof. 
\end{proof}
\vspace{5pt}

\begin{figure}[h]
\begin{center}
\includegraphics[scale=.4]{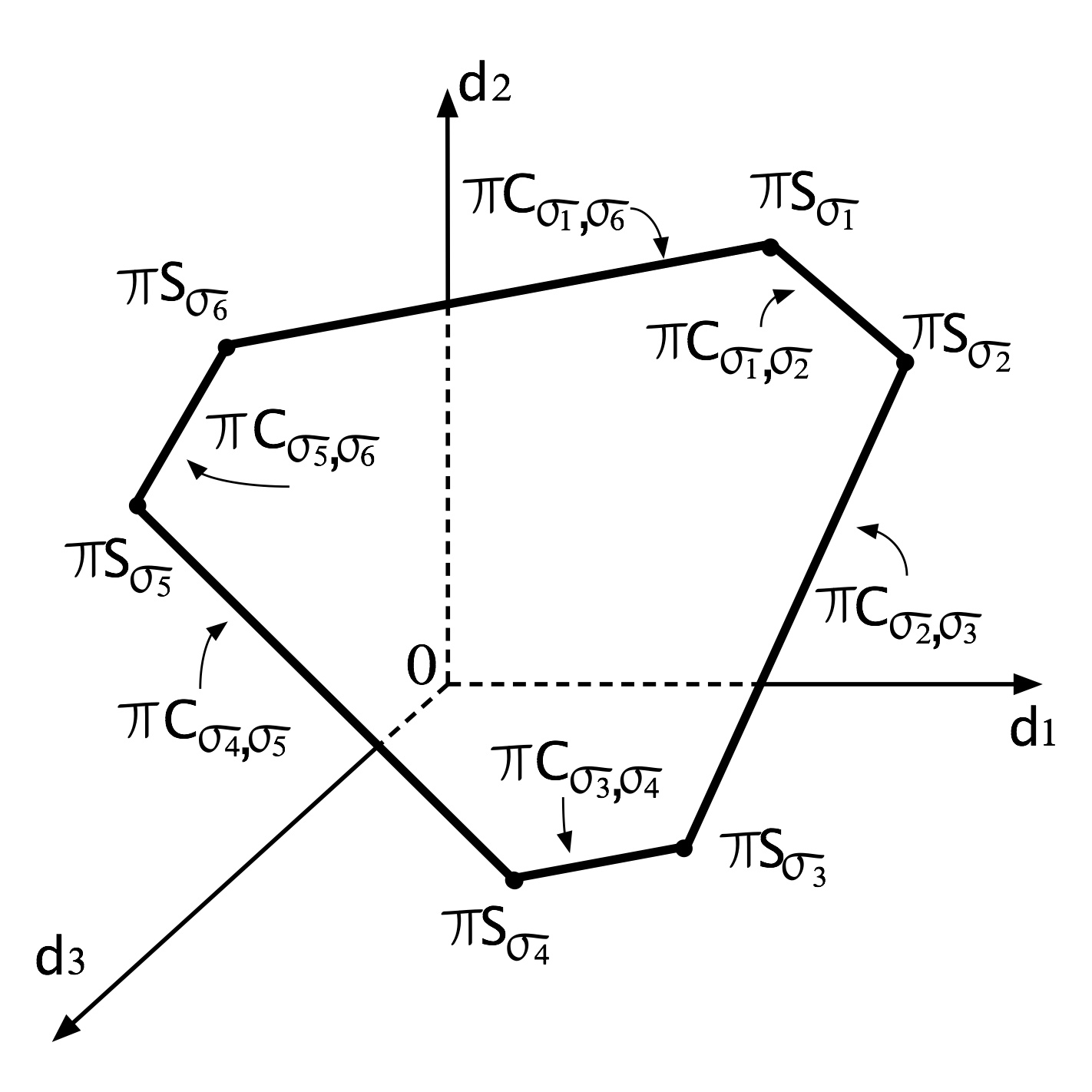}
\caption {Each vertex is the projection of a diagonal matrix and each edge is the projection of an unstable manifold $W^u(K_{\sigma_i,\sigma_j})$. For convenience, we use $C_{\sigma_i,\sigma_j}$ to denote $W^u(K_{\sigma_i,\sigma_j})$ with emphasis that each $W^u(\sigma_i,\sigma_j)$ consists of two disjoint curves with the same image under $\pi$.}
\label{FIG3}
\end{center}
\end{figure}

There are as many as $\frac{1}{2}(n-1)\cdot n!$ pairs of adjacent neighbors. To better understand the geometry behind the theorem, we consider the convex polytope $C(\Lambda)$. Each vertex of the polytope, known as a $0$-dim face, corresponds to a vector $(\lambda_{\sigma(1)},\cdots,\lambda_{\sigma(n)})$. It is the image of a diagonal matrix under the projection map $\pi$. Each edge of the polytope, known as a $1$-dim face, is then the image of an unstable manifold $W^u(K_{\sigma,\sigma'})$. Notice that each $W^u(K_{\sigma,\sigma'})$ has two disjoint curves, and both have the same image under $\pi$. It is clear then each edge corresponds to a pair of adjacent neighbors.

\section{Application of adjacency theorem on stochastic gradient flow}
There is a stochastic version of the double bracket flow by adding an isotropic noise into the equation
\begin{align}
dH= & [H,[H,\pi(H)]]dt + \epsilon \sum[\Omega_{ij},H]d\omega_{ij}\notag\\
& +\frac{\epsilon^2}{2}\sum[\Omega_{ij},[\Omega_{ij},H]]dt
\end{align}
The third term in the right hand side of equation above appears as a consequence of the It\^o rule so that  the solution still evolves on the isospectral manifold. In the literatures of NMR, the last two terms relate to the Lindblad terms modeling the heat bath. Each $\Omega_{ij}$ is a skew symmetric matrix defined by $\Omega_{ij}=\vec e_i\vec e^T_j-\vec e_j\vec e^T_i$ where $\{\vec e_1,\cdots,\vec e_n\}$ is a standard basis in $\mathbb{R}^n$. The significance of the stochastic effects is modeled by the scalar $\epsilon$. 

It happens that there is an explicit formula for the steady state solution of the Fokker-Planck equation associated with our stochastic equation. It takes the form
\begin{equation}
\rho(H)=c\exp(\frac{2\Psi(H)}{\epsilon^2})
\end{equation}
where $c$ is a normalization factor of the density. (see, for example, \cite{brockett1997notes} the derivation of the formula.) This implies that the equiprobability surface of the stochastic flow coincides with the equipotential surface of $\Psi$. If $\epsilon$ is small enough, then the density function is highly peaked at diagonal matrices as they are the local maxima. So a typical trajectory will spend most of its time around stable equilibria, this suggests that we approximate the behavior of the sample path by setting up a Markov model whose states are the $n!$ diagonal matrices, and the trajectory of a sample path is simplified by a chain 
\begin{equation}
s_{\sigma_1}\xrightarrow{T_1}s_{\sigma_2}\xrightarrow{T_2}\cdots\xrightarrow{T_{k-1}}s_{\sigma_k}
\end{equation} 
The rest of this section is to develop a method to evaluate the transition probabilities. 

There are two main challenging problems when coming to the computation: the scale and the model. For each state, there are as many as $(n!-1)$ transition probabilities we need to evaluate, so in general the amount of computation is about $n!(n!-1)$. To evaluate each transition probability, we need to investigate a corresponding first hitting time model: we ask for the probability $P(T|s_{\sigma_i}\to s_{\sigma_j})$ of the event that the passage time is than $T$ for a sample path to reach state $s_{\sigma_j}$ from state $s_{\sigma_i}$ without visiting any other state. We now show how adjacency theorem enters to simplify the problems  
\vspace{5pt}

\subsubsection*{Reduce the computational complexity} Let $A_{\sigma_i}$ be the set of adjacent neighbors of $s_{\sigma_i}$. Suppose $s_{\sigma_j}\notin A_{\sigma_i}$, then we infer that under the case $\epsilon\ll 1$
\begin{equation}\label{PROB}
P(T|s_{\sigma_i}\to s_{\sigma_j})\ll \sum_{\sigma\in A_{\sigma_i}}P(T|s_{\sigma_i}\to s_{\sigma})
\end{equation}
for any $T>0$. This is reasonable because from geometric perspective, if a sample path escapes from the region of attraction $W^s(s_{\sigma_i})$, there is a high percentage a path gets trapped by one of its neighbors. So our first  step of approximation is to set $P(T|s_{\sigma_i}\to s_{\sigma_j})=0$ for any $\sigma_j\notin A_{\sigma_i}$. So then for each state, the amount of computation is reduced to $(n-1)$.

\begin{figure}[h]
\begin{center}
\includegraphics[scale=.4]{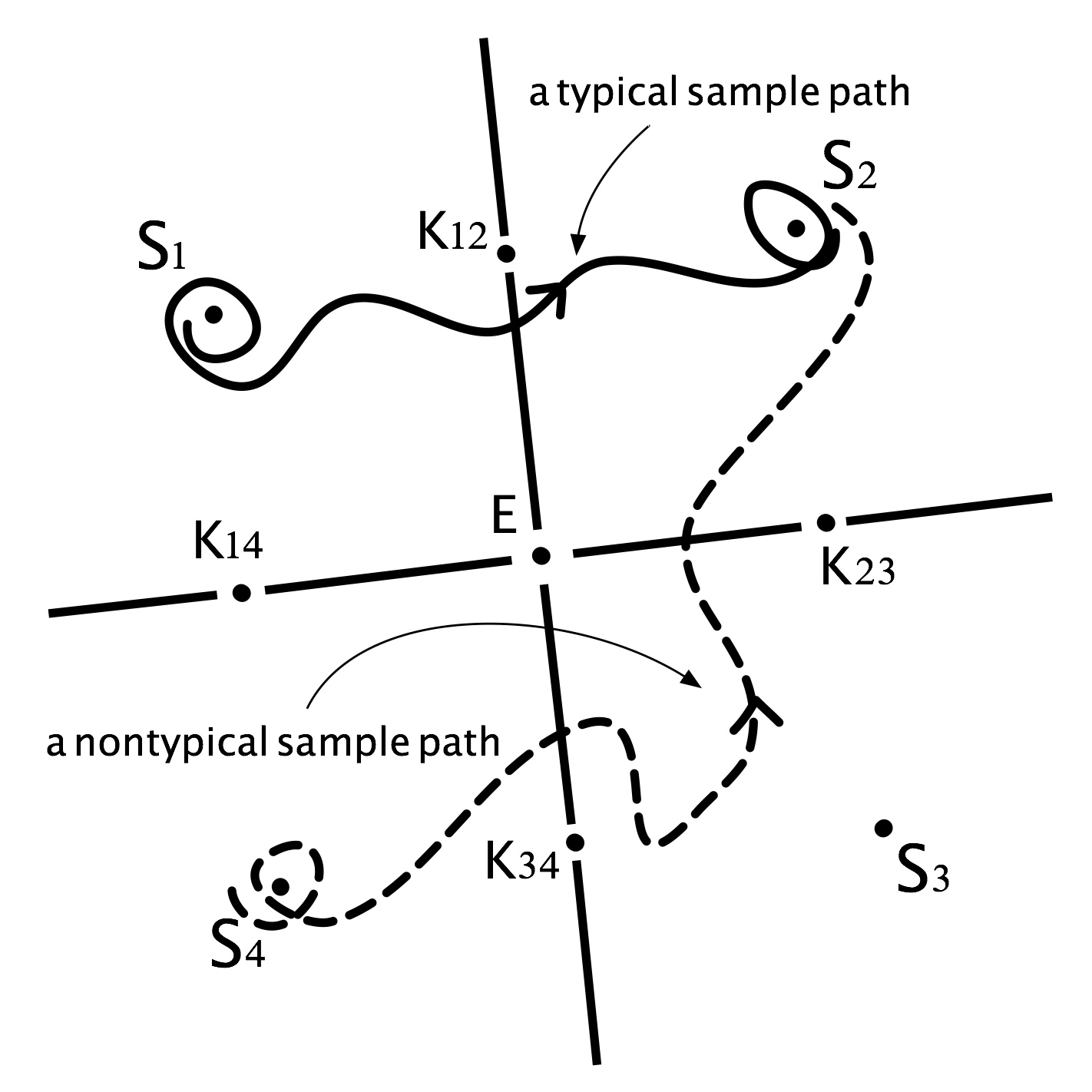}
\caption{Illustrating the idea behind equation \eqref{PROB}: a typical sample path only connects adjacent neighbors while a nontypical sample path does not.}
\end{center}
\end{figure}

\subsubsection*{Approximate the transition probability} For each state $s_{\sigma_i}$, the sample space of transition probability is discrete in space $A_{\sigma_i}$ but continuous in time. In general there is no exact formula for computing $P(T|s_{\sigma_i}\to s_{\sigma_j})$. One approach to approximate the density function is to relate the first hitting problem to an optimal control problem. 

In quantum mechanics, it is known that the probability for a system to stay in a quantum state of energy $E$ is proportional to $\exp(-E)$. This idea of Boltzmann sheds light on our problem.  We consider the control problem 
\begin{equation}
\dot H=[H,[H,\pi(H)]]+\epsilon[H,U]
\end{equation}
with $U$ skew symmetric and the goal is to minimize the energy, i.e,
\begin{equation}
E(T|\sigma_i\to\sigma_j):=-\min_{U(t)}\frac{1}{2}\int^T_{0}tr(U^2(t))dt
\end{equation}
under the assumption that $H(0)=\sigma_{i}$ and $H(T)=\sigma_j$ are fixed. We then approximate $P(T|\sigma_{i}\to\sigma_{j})$ by a scalar proportion of $\exp(-E(T|\sigma_i\to\sigma_j))$.

By solving the Euler-Lagrange equation, we conclude that each optimal trajectory, or the \textbf{energy minimizing path}(EMP) has to satisfy
\begin{align}
& \dot{H}=[H,\Omega] \\
& \dot{\Omega}=[H,[\pi(H),[H,\pi(H)]]]+[H,\pi([H,[H,\pi(H)]])]
\end{align}
It is hard to compute the MEP in general, however in the case $H(0)$ and $H(T)$ are both diagonal matrices, each MEP coincides with a geodesic. In particular, there are two MEPs and they together form the unstable manifold of  $K_{\sigma_1,\sigma_2}$. This then simplifies the situation to a scalar problem: suppose the simple transposition relating $\sigma_i$ and $\sigma_j$ is the 2-cycle $(\lambda_1,\lambda_2)$, then the control model is given by
\begin{equation}
\dot\theta=-\frac{1}{2}(\lambda_1-\lambda_2)^2\sin(2\theta)+2\epsilon u
\end{equation} 
and the goal is to minimize $\int^T_0u^2dt$. 
\vspace{5pt}  

 Before we ending this section, we point out that locating a MEP that connects two adjacent neighbors and computing the minimal consumption of energy is more than an ad-hoc plan for evaluating the transition probability. For example, consider the situation that the control is intermittent and impulse-like, and our goal is to steer the system from one stable equilibrium to the other. Then questions such as in what direction one can escape from a region of attraction by means of an impulse? How we can save most of energy and/or time to reach the target? These questions relate to the design of a path that concatenates pairs of adjacent neighbors, the analysis done in this section is then essential. 

\section*{Acknowledgements}

The author thanks Dr. Roger W. Brockett at Harvard University for his comments on an earlier draft.


\bibliographystyle{unsrt}
\bibliography{adjacencycriterion}

\end{document}